\documentclass{article}%
\usepackage{amssymb}
\usepackage{graphicx}
\usepackage{amsmath}%
\usepackage{amsfonts}
\providecommand{\U}[1]{\protect\rule{.1in}{.1in}}

\begin{document}

\title{A Spectral Method for Elliptic Equations: \\The Neumann Problem}
\author{Kendall Atkinson\\Departments of Mathematics \& Computer Science \\The University of Iowa
\and David Chien, Olaf Hansen\\Department of Mathematics \\California State University San Marcos}
\maketitle

\begin{abstract}
Let $\Omega$ be an open, simply connected, and bounded region in
$\mathbb{R}^{d}$, $d\geq2$, and assume its boundary $\partial\Omega$ is
smooth. Consider solving an elliptic partial differential equation $-\Delta
u+\gamma u=f$ over $\Omega$ with a Neumann boundary condition. The problem is
converted to an equivalent\ elliptic problem over the unit ball $B$, and then
a spectral Galerkin method is used to create a convergent sequence of
multivariate polynomials $u_{n}$ of degree $\leq n$ that is convergent to $u$.
The transformation from $\Omega$ to $B$\ requires a special analytical
calculation for its implementation. With sufficiently smooth problem
parameters, the method is shown to be rapidly convergent. For $u\in C^{\infty
}\left(  \overline{\Omega}\right)  $ and assuming $\partial\Omega$ is a
$C^{\infty}$ boundary, the convergence of $\left\Vert u-u_{n}\right\Vert
_{H^{1}}$ \ to zero is faster than any power of $1/n$. Numerical examples in
$\mathbb{R}^{2}$ and $\mathbb{R}^{3}$ show experimentally an exponential rate
of convergence.

\end{abstract}

\section{INTRODUCTION}

Consider\ solving the Neumann problem for Poisson's equation:
\begin{align}
-\Delta u+\gamma(s)u &  =f(s),\quad\quad s\in\Omega\medskip\label{en1}\\
\frac{\partial u(s)}{\partial n_{s}} &  =g(s),\quad\quad s\in\partial
\Omega.\label{en2}%
\end{align}
Assume $\Omega$ is an open, simply-connected, and bounded region in
$\mathbb{R}^{d}$, $d\geq2$, and assume that its boundary $\partial\Omega$ is
several times continuously differentiable. Similarly, assume the functions
$\gamma(s)$ and $f(s)$ are several times continuously differentiable over
$\overline{\Omega}$, and assume that $g(s)$ is several times continuously
differentiable over the boundary $\partial\Omega$.

There is a rich literature on spectral methods for solving partial
differential equations. From the more recent literature, we cite \cite{cqhz1},
\cite{cqhz2}, \cite{doha}, and \cite{shen}. Their bibliographies contain
\ references to earlier papers on spectral methods. The present paper is a
continuation of the work in \cite{ACH-Dir} in which a spectral method is given
for a general elliptic equation with a Dirichlet boundary condition. Our
approach is somewhat different than the standard approaches. We convert the
partial differential equation to an equivalent problem on the unit disk or
unit ball, and in the process we are required to work with a more complicated
equation. Our approach is reminiscent of the use of conformal mappings for
planar problems. Conformal mappings can be used with our approach when working
on planar problems, although having a conformal mapping is not necessary.

In \S \ref{unique} we assume that (\ref{en1})-(\ref{en2}) is uniquely
solvable, and we present a spectral Galerkin method for its solution. In
\S \ref{nonunique} we extend the method to the problem with $\gamma(s)\equiv0$
in $\Omega$. The problem is no longer uniquely solvable and we extend our
spectral method to this case. The implementation of the method is discussed in
\S \ref{implementation} and it is illustrated in \S \ref{examples}.

\section{A spectral method for the uniquely solvable case \label{unique}}

We assume the Neumann problem (\ref{en1})-(\ref{en2}) is uniquely solvable.
This is true, for example, if
\begin{equation}
\gamma\left(  s\right)  \geq c_{\gamma}>0,\quad\quad s\in\overline{\Omega
}\label{en2a}%
\end{equation}
for some constant $c_{\gamma}>0$. For functions $u\in H^{2}\left(
\Omega\right)  ,$ $v\in H^{1}\left(  \Omega\right)  $,%
\begin{equation}%
\begin{array}
[c]{r}%
{\displaystyle\int_{\Omega}}
v(s)\left[  -\Delta u(s)+\gamma(s)u\right]  \,ds=%
{\displaystyle\int_{\Omega}}
\left[  \triangledown u(s)\cdot\triangledown v(s)+\gamma(s)u(s)v(s)\right]
\,ds\quad\medskip\\
-%
{\displaystyle\int_{\partial\Omega}}
v\left(  s\right)  \dfrac{\partial u(s)}{\partial n_{s}}\,ds.
\end{array}
\label{en3}%
\end{equation}
Introduce the bilinear functional%
\begin{equation}
\mathcal{A}\left(  v_{1},v_{2}\right)  =\int_{\Omega}\left[  \triangledown
v_{1}(s)\cdot\triangledown v_{2}(s)+\gamma(s)v_{1}(s)v_{2}(s)\right]
\,ds.\label{en4}%
\end{equation}
The variational form of the Neumann problem (\ref{en1})-(\ref{en2}) is as
follows: find $u$ such that%
\begin{equation}
\mathcal{A}\left(  u,v\right)  =\ell_{1}(v)+\ell_{2}\left(  v\right)
,\quad\quad\forall v\in H^{1}\left(  \Omega\right) \label{en8}%
\end{equation}
with the linear functionals defined by%
\begin{align}
\ell_{1}(v) &  =\int_{\Omega}v(s)f(s)\,ds,\medskip\label{en8a}\\
\ell_{2}\left(  v\right)   &  =\int_{\partial\Omega}v\left(  s\right)
g(s)\,ds.\label{en8b}%
\end{align}
The norms we use for $\ell_{1}$ and $\ell_{2}$ are the standard operator norms
when regarding $\ell_{1}$ and $\ell_{2}$ as linear functionals on
$H^{1}\left(  \Omega\right)  $. The functional $\ell_{1}$ is bounded easily on
$H^{1}\left(  \Omega\right)  $,%
\begin{equation}
\left\vert \ell_{1}\left(  v\right)  \right\vert \leq\Vert f\Vert_{L^{2}}\Vert
v\Vert_{L^{2}}\leq\Vert f\Vert_{L^{2}}\Vert v\Vert_{H^{1}}.\label{en9}%
\end{equation}
Ordinarily, we will use $\Vert v\Vert_{1}$ in place of $\Vert v\Vert_{H^{1}}$.

The functional $\ell_{2}$ is
bounded (at least for bounded domains $\Omega$). To show this, begin by noting
that the restriction $\rho:H^{1}(\Omega)\rightarrow H^{1/2}(\partial\Omega)$
is continuous \cite[Th. 3.37]{McLean} and the imbedding $\iota:H^{1/2}%
(\partial\Omega)\hookrightarrow L^{2}(\partial\Omega) $ is compact \cite[Th.
3.27]{McLean}. If we further denote by $l_{g}$ the continuous mapping
\[
l_{g}:u\mapsto\int_{\partial\Omega}u(s)g(s)\,ds,\quad\quad u\in L^{2}%
(\partial\Omega)
\]
then we see $\ell_{2}=l_{g}\circ\iota\circ\rho$, and therefore $\ell_{2}$ is
bounded.

It is straightforward to show $\mathcal{A}$ is bounded,%
\begin{equation}%
\begin{array}
[c]{c}%
\left\vert \mathcal{A}\left(  v,w\right)  \right\vert \leq c_{\mathcal{A}%
}\left\Vert v\right\Vert _{1}\left\Vert w\right\Vert _{1},\medskip\\
c_{\mathcal{A}}=\max\left\{  1,\left\Vert \gamma\right\Vert _{\infty}\right\}
.
\end{array}
\label{en13}%
\end{equation}
In addition, we assume $\mathcal{A}$ is strongly elliptic on $H^{1}\left(
\Omega\right)  $,
\begin{equation}
\mathcal{A}\left(  v,v\right)  \geq c_{e}\Vert v\Vert_{1}^{2},\quad\quad v\in
H^{1}\left(  \Omega\right) \label{en23}%
\end{equation}
with some $c_{e}>0$. This follows ordinarily from showing the unique
solvability of the Neumann problem (\ref{en1})-(\ref{en2}). If (\ref{en2a}) is
satisfied, then we can satisfy (\ref{en23}) with%
\[
c_{e}=\min\left\{  1,c_{\gamma}\right\}
\]
Under our assumptions on $\mathcal{A}$, including the strong ellipticity in
(\ref{en23}), the Lax-Milgram Theorem implies the existence of a unique
solution $u$ to (\ref{en8}) with
\begin{equation}
\Vert u\Vert_{1}\leq\frac{1}{c_{e}}\left[  \Vert\ell_{1}\Vert+\Vert\ell
_{2}\Vert\right]  .\label{en15}%
\end{equation}

Our spectral method is defined using polynomial approximations over the open
unit ball in $\mathbb{R}^{d}$, call it $B_{d}$. Introduce a change of
variables%
\[
\Phi:\overline{B}_{d}\underset{onto}{\overset{1-1}{\longrightarrow}}%
\overline{\Omega}%
\]
with $\Phi$ a twice-differentiable mapping, and let $\Psi=\Phi^{-1}%
:\overline{\Omega}\underset{onto}{\overset{1-1}{\longrightarrow}}\overline
{B}_{d}$. [We comment later on the creation of $\Phi$ for cases in which only
the boundary mapping $\phi:\partial B_{d}\rightarrow\partial\Omega$ is known.]
\ For $v\in L^{2}\left(  \Omega\right)  $, let%
\[
\widetilde{v}(x)=v\left(  \Phi\left(  x\right)  \right)  ,\quad\quad
x\in\overline{B}_{d}\subseteq\mathbb{R}^{d}%
\]
and conversely,%
\[
v(s)=\widetilde{v}\left(  \Psi\left(  s\right)  \right)  ,\quad\quad
s\in\overline{\Omega}\subseteq\mathbb{R}^{d}.
\]
Assuming $v\in H^{1}\left(  \Omega\right)  $, we can show%
\[
\nabla_{x}\widetilde{v}\left(  x\right)  =J\left(  x\right)  ^{\text{T}}%
\nabla_{s}v\left(  s\right)  ,\quad\quad s=\Phi\left(  x\right)
\]
with $J\left(  x\right)  $ the Jacobian matrix for $\Phi$ over the unit ball
$B_{d}$,%
\[
J(x)\equiv\left(  D\Phi\right)  (x)=\left[  \frac{\partial\Phi_{i}%
(x)}{\partial x_{j}}\right]  _{i,j=1}^{d},\quad\quad x\in\overline{B}_{d}.
\]
Similarly,%
\[
\nabla_{s}v(s)=K(s)^{\text{T}}\nabla_{x}\widetilde{v}(x),\quad\quad x=\Psi(s)
\]
with $K(s)$ the Jacobian matrix for $\Psi$ over $\Omega$. Also,
\begin{equation}
K\left(  \Phi\left(  x\right)  \right)  =J\left(  x\right)  ^{-1}.\label{en16}%
\end{equation}

Using the change of variables $s=\Phi\left(  x\right)  $, the formula
(\ref{en4}) converts to
\begin{align}
\mathcal{A}\left(  v_{1},v_{2}\right)   &  =\int_{B_{d}}\{[K\left(
\Phi\left(  x\right)  \right)  ^{\text{T}}\nabla_{x}\widetilde{v}_{1}\left(
x\right)  ]^{\text{T}}[K\left(  \Phi\left(  x\right)  \right)  ^{\text{T}%
}\nabla_{x}\widetilde{v}_{2}\left(  x\right)  ]\smallskip\nonumber\\
&  \left.  \quad\quad\right.  +\gamma(\Phi\left(  x\right)  )v_{1}(\Phi\left(
x\right)  )v_{2}(\Phi\left(  x\right)  \}\,\left\vert \det\left[  J(x)\right]
\right\vert \,dx\smallskip\nonumber\\
&  =\int_{B_{d}}\{[J\left(  x\right)  ^{-\text{T}}\nabla_{x}\widetilde{v}%
_{1}\left(  x\right)  ]^{\text{T}}[J\left(  x\right)  ^{-\text{T}}\nabla
_{x}\widetilde{v}_{2}\left(  x\right)  ]\smallskip\nonumber\\
&  \left.  \quad\quad\right.  +\widetilde{\gamma}(x)\widetilde{v}%
_{1}(x)\widetilde{v}_{2}(x)\}\,\left\vert \det\left[  J(x)\right]  \right\vert
\,dx\smallskip\nonumber\\
&  =\int_{B_{d}}\{\nabla_{x}\widetilde{v}_{1}\left(  x\right)  ^{\text{T}%
}A(x)\nabla_{x}\widetilde{v}_{2}\left(  x\right)  +\widetilde{\gamma
}(x)\widetilde{v}_{1}(x)\widetilde{v}_{2}(x)\}\left\vert \det\left[
J(x)\right]  \right\vert \,dx\smallskip\nonumber\\
&  \equiv\widetilde{\mathcal{A}}\left(  \widetilde{v}_{1},\widetilde{v}%
_{2}\right) \label{en20}%
\end{align}
with%
\[
A(x)=J\left(  x\right)  ^{-1}J\left(  x\right)  ^{-\text{T}}.
\]

We can also introduce analogues to $\ell_{1}$ and $\ell_{2}$ following a
change of variables, calling them $\widetilde{\ell}_{1}$ and $\widetilde{\ell
}_{2}$ and defined on $H^{1}\left(  B_{d}\right)  $. For example,%
\[
\widetilde{\ell}_{1}(\widetilde{v})=\int_{B_{d}}\widetilde{v}(x)f(\Phi
(x))\left\vert \det\left[  J(x)\right]  \right\vert \,dx.
\]
We can then convert (\ref{en8}) to an equivalent problem over $H^{1}\left(
B_{d}\right)  $. The variational problem becomes
\begin{equation}
\widetilde{\mathcal{A}}\left(  \widetilde{u},\widetilde{v}\right)
=\widetilde{\ell}_{1}(\widetilde{v})+\widetilde{\ell}_{2}\left(  \widetilde
{v}\right)  ,\quad\quad\forall\widetilde{v}\in H^{1}\left(  B_{d}\right)
.\label{en21}%
\end{equation}
The assumptions and results in (\ref{en8})-(\ref{en23}) extend to this new
problem on $H^{1}\left(  B_{d}\right)  $. The strong ellipticity condition
(\ref{en23}) becomes%
\begin{align}
\widetilde{\mathcal{A}}\left(  \widetilde{v},\widetilde{v}\right)   &
\geq\widetilde{c}_{e}\Vert\widetilde{v}\Vert_{1}^{2},\quad\quad\widetilde
{v}\in H^{1}\left(  B_{d}\right)  ,\medskip\label{en22}\\
\widetilde{c}_{e} &  =c_{e}\frac{\min_{x\in\overline{B}_{d}}\left\vert \det
J(x)\right\vert }{\max\left[  1,\max_{x\in\overline{B}_{d}}\left\Vert
J(x)\right\Vert _{2}^{2}\right]  }\nonumber
\end{align}
where $\left\Vert J(x)\right\Vert _{2}$ denotes the operator matrix 2-norm of
$J(x)$ for $\mathbb{R}^{d}$. Also,%
\[%
\begin{array}
[c]{c}%
\left\vert \widetilde{\mathcal{A}}\left(  \widetilde{v},\widetilde{w}\right)
\right\vert \leq\widetilde{c}_{\mathcal{A}}\left\Vert \widetilde{v}\right\Vert
_{1}\left\Vert \widetilde{w}\right\Vert _{1},\medskip\\
\widetilde{c}_{\mathcal{A}}=\left\{  \max\limits_{x\in\overline{B}_{d}%
}\left\vert \det\left[  J(x)\right]  \right\vert \right\}  \,\max\left\{
\max\limits_{x\in\overline{B}_{d}}\left\Vert A(x)\right\Vert _{2},\left\Vert
\gamma\right\Vert _{\infty}\right\}  .
\end{array}
\]

For the finite dimensional problem, we want to use the approximating subspace
$\Pi_{n}\equiv\Pi_{n}^{d}$. We want to find $\widetilde{u}_{n}\in\Pi_{n}$ such
that%
\begin{equation}
\widetilde{\mathcal{A}}\left(  \widetilde{u}_{n},\widetilde{v}\right)
=\widetilde{\ell}_{1}(\widetilde{v})+\widetilde{\ell}_{2}\left(  \widetilde
{v}\right)  ,\quad\quad\forall\widetilde{v}\in\Pi_{n}.\label{en25}%
\end{equation}
The Lax-Milgram Theorem (cf. \cite[\S 8.3]{atkinson-han}, \cite[\S 2.7]%
{BrennerScott}) implies the existence of $u_{n}$ for all $n$. For the error in
this Galerkin method, Cea's Lemma (cf. \cite[p. 365]{atkinson-han}, \cite[p.
62]{BrennerScott}) implies the convergence of $u_{n}$ to $u$, and moreover,%
\begin{equation}
\Vert\widetilde{u}-\widetilde{u}_{n}\Vert_{1}\leq\frac{\widetilde
{c}_{\mathcal{A}}}{\widetilde{c}_{e}}\inf_{\widetilde{v}\in\Pi_{n}}%
\Vert\widetilde{u}-\widetilde{v}\Vert_{1}.\label{en31}%
\end{equation}
It remains to bound the best approximation error on the right side of this inequality.

Ragozin \cite{ragozin} gives bounds on the rate of convergence of best
polynomial approximation over the unit ball, and these results are extended in
\cite{bbl} to simultaneous approximation of a function and some of its lower
order derivatives. Assume $\widetilde{u}\in C^{m+1}\left(  \overline{B}%
_{d}\right)  $. Using \cite[Theorem 1]{bbl}, we have
\begin{equation}
\inf_{\widetilde{v}\in\Pi_{n}}\Vert\widetilde{u}-\widetilde{v}\Vert_{1}%
\leq\frac{c(u,m)}{n^{m}}\omega_{u,m+1}\left(  \frac{1}{n}\right) \label{en33}%
\end{equation}
with%
\[
\omega_{u,m+1}\left(  \delta\right)  =\sup_{\left\vert \alpha\right\vert
=m+1}\left(  \sup_{\left\vert x-y\right\vert \leq\delta}\left\vert D^{\alpha
}\widetilde{u}\left(  x\right)  -D^{\alpha}\widetilde{u}\left(  y\right)
\right\vert \right)  .
\]
The notation $D^{\alpha}\widetilde{u}\left(  x\right)  $ is standard
derivative notation with $\alpha$ a multi-integer. \ In particular, for
$\alpha=\left(  \alpha_{1},\dots,\alpha_{d}\right)  $,%
\[
D^{\alpha}\widetilde{u}\left(  x\right)  =\frac{\partial^{\left\vert
\alpha\right\vert }\widetilde{u}\left(  x_{1},\dots,x_{d}\right)  }{\partial
x_{1}^{\alpha_{1}}\cdots\partial x_{d}^{\alpha_{d}}}.
\]
When (\ref{en33}) is combined with (\ref{en31}), we see that our solutions
$\widetilde{u}_{n}$ converge faster than any power of $1/n$ provided
$\widetilde{u}\in C^{\infty}\left(  \overline{B}_{d}\right)  $.

\section{A spectral method for $-\Delta u=f$ \label{nonunique}}

Consider the Neumann problem for Poisson's equation:
\begin{align}
-\Delta u  &  =f(s),\quad\quad s\in\Omega\medskip\label{en40}\\
\frac{\partial u(s)}{\partial n_{s}}  &  =g(s),\quad\quad s\in\partial
\Omega\label{en41}%
\end{align}
As a reference for this problem, see \cite[\S 5.2]{BrennerScott}.

As earlier in (\ref{en3}), we have for functions $u\in H^{2}\left(
\Omega\right)  ,$ $v\in H^{1}\left(  \Omega\right)  $,
\begin{equation}
\int_{\Omega}v(s)\Delta u(s)\,ds=-\int_{\Omega}\triangledown u(s)\cdot
\triangledown v(s)\,ds+\int_{\partial\Omega}v\left(  s\right)  \frac{\partial
u(s)}{\partial n_{s}}\,ds\label{en42}%
\end{equation}
If this Neumann problem (\ref{en40})-(\ref{en41}) is solvable, then its
solution is not unique: any constant added to a solution gives another
solution. In addition, if (\ref{en40})-(\ref{en41}) is solvable, then
\begin{equation}
\int_{\Omega}v(s)f(s)\,ds=\int_{\Omega}\triangledown u(s)\cdot\triangledown
v(s)\,ds-\int_{\partial\Omega}v\left(  s\right)  g(s)\,ds\label{en43}%
\end{equation}
Choosing $v(s)\equiv1$, we obtain%
\begin{equation}
\int_{\Omega}f(s)\,ds=-\int_{\partial\Omega}g(s)\,ds\label{en5}%
\end{equation}
This is a necessary and sufficient condition on the functions $f$ and $g$ in
order that (\ref{en40})-(\ref{en41}) be solvable. With this constraint, the
Neumann problem is solvable. To deal with the non-unique solvability, we look
for a solution $u$ satisfying%
\begin{equation}
\int_{\Omega}u(s)\,ds=0\label{en6}%
\end{equation}

Introduce the bilinear functional%
\begin{equation}
\mathcal{A}\left(  v_{1},v_{2}\right)  =\int_{\Omega}\triangledown
v_{1}(s)\cdot\triangledown v_{2}(s)\,ds\label{en6.a}%
\end{equation}
and the function space%
\begin{equation}
\mathcal{V}=\left\{  v\in H^{1}\left(  \Omega\right)  :\int_{\Omega
}v(s)\,ds=0\right\} \label{en7}%
\end{equation}
$\mathcal{A}$ is bounded,%
\[
\left\vert \mathcal{A}\left(  v,w\right)  \right\vert \leq\left\Vert
v\right\Vert _{1}\left\Vert w\right\Vert _{1},\quad\quad\forall v,w\in
\mathcal{V}.
\]
From \cite[Prop. 5.3.2]{BrennerScott} $\mathcal{A}\left(  \cdot,\cdot\right)
$ is strongly elliptic on $\mathcal{V}$, satisfying
\[
\mathcal{A}\left(  v,v\right)  \geq c_{e}\Vert v\Vert_{1}^{2},\quad\quad
v\in\mathcal{V}%
\]
for some $c_{e}>0$. The variational form of the Neumann problem (\ref{en40}%
)-(\ref{en41}) is as follows: find $u$ such that%
\begin{equation}
\mathcal{A}\left(  u,v\right)  =\ell_{1}(v)+\ell_{2}\left(  v\right)
,\quad\quad\forall v\in\mathcal{V}\label{en45}%
\end{equation}
with $\ell_{1}$ and $\ell_{2}$ defined as in (\ref{en8a})-(\ref{en8b}). As
before, the Lax-Milgram Theorem implies the existence of a unique solution $u
$ to (\ref{en45}) with
\[
\Vert u\Vert_{1}\leq\frac{1}{c_{e}}\left[  \Vert\ell_{1}\Vert+\Vert\ell
_{2}\Vert\right]  .
\]

As in the preceding section, we transform the problem from being defined over
$\Omega$ to being over $B_{d}$. Most of the arguments are repeated, and we
have%
\[
\widetilde{\mathcal{A}}\left(  \widetilde{v}_{1},\widetilde{v}_{2}\right)
=\int_{B_{d}}\{\nabla_{x}\widetilde{v}_{1}\left(  x\right)  ^{\text{T}%
}A(x)\nabla_{x}\widetilde{v}_{2}\left(  x\right)  \}\left\vert \det\left[
J(x)\right]  \right\vert \,dx.
\]
The condition (\ref{en6}) becomes%
\[
\int_{B}\widetilde{v}(x)\left\vert \det\left[  J(x)\right]  \right\vert
\,dx=0.
\]
We introduce the space
\begin{equation}
\widetilde{\mathcal{V}}=\left\{  \widetilde{v}\in H^{1}\left(  B\right)
:\int_{B}\widetilde{v}(x)\left\vert \det\left[  J(x)\right]  \right\vert
\,dx=0\right\}  .\label{en48}%
\end{equation}
The Neumann problem now has the reformulation%
\begin{equation}
\widetilde{\mathcal{A}}\left(  \widetilde{u},\widetilde{v}\right)
=\widetilde{\ell}_{1}(\widetilde{v})+\widetilde{\ell}_{2}\left(  \widetilde
{v}\right)  ,\quad\quad\forall\widetilde{v}\in\widetilde{\mathcal{V}%
}\label{en49}%
\end{equation}

For the finite dimensional approximating problem, we use
\begin{equation}
\widetilde{\mathcal{V}}_{n}=\widetilde{\mathcal{V}}\cap\Pi_{n}\label{en51}%
\end{equation}
Then we want to find $\widetilde{u}_{n}\in\widetilde{\mathcal{V}}_{n}$ such
that%
\begin{equation}
\widetilde{a}\left(  \widetilde{u}_{n},\widetilde{v}\right)  =\widetilde{\ell
}_{1}(\widetilde{v})+\widetilde{\ell}_{2}\left(  \widetilde{v}\right)
,\quad\quad\forall\widetilde{v}\in\widetilde{\mathcal{V}}_{n}\label{en52}%
\end{equation}
We can invoke the standard results of the Lax-Milgram Theorem and Cea's Lemma
to obtain the existence of a unique solution $\widetilde{u}_{n}$, and
moreover,%
\begin{equation}
\Vert\widetilde{u}-\widetilde{u}_{n}\Vert_{1}\leq c\inf_{v\in\widetilde
{\mathcal{V}}_{n}}\Vert\widetilde{u}-v\Vert_{1}.\label{en60}%
\end{equation}
for some $c>0$. A modification of the argument that led to (\ref{en33}) can be
used to obtained a similar result for (\ref{en60}). First, however, we discuss
the practical problem of choosing a basis for $\widetilde{\mathcal{V}}_{n}$.

\subsection{Constructing a basis for\ $\widetilde{\mathcal{V}}_{n}$}

Let $\left\{  \varphi_{j}:1\leq j\leq N_{n}^{d}\right\}  $ denote a basis for
$\Pi_{n}$ (usually we choose $\left\{  \varphi_{j}\right\}  $ to be an
orthogonal family in the norm of $L^{2}\left(  B_{d}\right)  $). \ We assume
that $\varphi_{1}(x)\ $is a nonzero constant function. Introduce the new basis
elements%
\begin{equation}
\widehat{\varphi}_{j}=\varphi_{j}-\frac{1}{C}\int_{B}\varphi_{j}(x)\left\vert
\det\left[  J(x)\right]  \right\vert \,dx,\quad\quad1\leq j\leq N_{n}%
^{d}\label{en80}%
\end{equation}
with
\begin{equation}
C=\int_{B}\left\vert \det\left[  J(x)\right]  \right\vert \,dx\equiv\left\Vert
\det\left[  J\right]  \right\Vert _{L^{1}}\label{en81}%
\end{equation}
Then $\widehat{\varphi}_{1}=0$ and
\begin{align*}
\int_{B}\widehat{\varphi}_{j}(x)\left\vert \det\left[  J(x)\right]
\right\vert \,dx &  =\int_{B}\varphi_{j}(x)\left\vert \det\left[  J(x)\right]
\right\vert \,dx\smallskip\\
&  -\frac{1}{C}\left[  \int_{B}\varphi_{j}(x)\left\vert \det\left[
J(x)\right]  \right\vert \,dx\right]  \left[  \int_{B}\left\vert \det\left[
J(x)\right]  \right\vert \,dx\right]  \smallskip\\
&  =0
\end{align*}
Thus $\left\{  \widehat{\varphi}_{j}:2\leq j\leq N_{n}^{d}\right\}  $ is a
basis of $\widetilde{\mathcal{V}}_{n}$ and we can use it for our Galerkin
procedure in (\ref{en52}).

\subsection{The rate of convergence of $\widetilde{u}_{n}$}

Now we estimate $\inf_{v\in\widetilde{\mathcal{V}}_{n}}\Vert\widetilde
{u}-v\Vert_{1}$; see (\ref{en60}). Recalling (\ref{en80}), we consider the
linear mapping $P:L^{2}(B_{d})\rightarrow L^{2}(B_{d})$ given by
\begin{align*}
(P\widetilde{u})(x) &  =\widetilde{u}(x)-\frac{1}{C}\int_{B}\left\vert
\det[J(y)]\right\vert \widetilde{u}(y)\,dy,\smallskip\\
C &  =\Vert\det[J]\Vert_{L^{1}};
\end{align*}
see (\ref{en81}). The mapping $P$ is a projection
\begin{align*}
P(P\widetilde{u})(x) &  =(P\widetilde{u})(x)-\frac{1}{C}\int_{B}\left\vert
\det[J(y)]\right\vert (P\widetilde{u})(y)\;dy\smallskip\\
&  =\widetilde{u}(x)-\frac{1}{C}\int_{B}\left\vert \det[J(y)]\right\vert
\;\widetilde{u}(y)\;dy-\\
&  \left(  \frac{1}{C}\int_{B}\left\vert \det[J(y)]\right\vert \left(
\widetilde{u}(y)-\frac{1}{C}\int_{B}\left\vert \det[J(z)]\right\vert
\;\widetilde{u}(z)\right)  \;dy\right)  \smallskip\\
&  =\widetilde{u}(x)-\frac{1}{C}\int_{B}\left\vert \det[J(y)]\right\vert
\;\widetilde{u}(y)\;dy-\frac{1}{C}\int_{B}\left\vert \det[J(y)]\right\vert
\widetilde{u}(y)\;dy\\
&  +\frac{1}{C^{2}}\underbrace{\int_{B}\left\vert \det[J(y)]\right\vert
\;dy}_{=C}\int_{B}\left\vert \det[J(z)]\right\vert \;\widetilde{u}%
(z)\;dz\smallskip\\
&  =\widetilde{u}(x)-\frac{1}{C}\int_{B}\left\vert \det[J(y)]\right\vert
\;\widetilde{u}(y)\;dy\smallskip\\
&  =(P\widetilde{u})(x)
\end{align*}
So $P^{2}=P$ and $P$ is a projection with $\Vert P\Vert_{L^{2}\rightarrow
L^{2}}\geq1$ and
\begin{align*}
\Vert P\widetilde{u}\Vert_{2} &  =\Vert\widetilde{u}-\frac{1}{C}\int
_{B}\left\vert \det[J(y)]\right\vert \,\widetilde{u}(y)\text{\ }dy\Vert
_{L^{2}}\smallskip\\
&  \leq\Vert\widetilde{u}\Vert_{L^{2}}+\frac{1}{C}\,\left\vert \int
_{B}\left\vert \det[J(y)]\right\vert \,\widetilde{u}(y)\,dy\right\vert
\,\Vert1\Vert_{L^{2}}\smallskip\\
&  \leq\Vert\widetilde{u}\Vert_{L^{2}}+\frac{1}{C}\Vert\det[J]\Vert_{L^{2}%
}\Vert\widetilde{u}\Vert_{L^{2}}\,\sqrt{\frac{\pi^{d/2}}{\Gamma\left(
1+\tfrac{1}{2}d\right)  }}\makebox[1cm]{}\mbox{(Cauchy-Schwarz)}\smallskip\\
&  =\left(  1+\sqrt{\frac{\pi^{d/2}}{\Gamma\left(  1+\tfrac{1}{2}d\right)  }%
}\frac{\Vert\det[J]\Vert_{L^{2}}}{\Vert\det[J]\Vert_{L^{1}}}\right)
\Vert\widetilde{u}\Vert_{L^{2}}\smallskip\\
&  =c_{P}\Vert\widetilde{u}\Vert_{L^{2}}%
\end{align*}
which shows $\Vert P\Vert_{L^{2}\rightarrow L^{2}}\leq c_{P}$ and
$\widetilde{\mathcal{V}}:=P(H^{1}(B_{d}))$, see (\ref{en48}). For
$\widetilde{u}\in H^{1}(B_{d})$ we also have $P\widetilde{u}\in H^{1}(B_{d})$
and here we again estimate the norm of $P$:
\begin{align*}
\Vert P\widetilde{u}\Vert_{H^{1}}^{2} &  =\Vert P\widetilde{u}\Vert_{2}%
^{2}+\Vert\nabla(P\widetilde{u})\Vert_{2}^{2}\smallskip\\
&  \leq c_{P}^{2}\Vert\widetilde{u}\Vert_{2}^{2}+\Vert\nabla\widetilde{u}%
\Vert_{2}^{2}%
\end{align*}
since $\nabla(P\widetilde{u})=\nabla\widetilde{u}$. Furthermore $c_{P}\geq1$,
so
\begin{align*}
\Vert P\widetilde{u}\Vert_{H^{1}}^{2} &  \leq c_{P}^{2}\Vert\widetilde{u}%
\Vert_{L^{2}}^{2}+c_{P}^{2}\Vert\nabla\widetilde{u}\Vert_{L^{2}}^{2}%
\smallskip\\
&  =c_{P}^{2}(\Vert\widetilde{u}\Vert_{L^{2}}^{2}+\Vert\nabla\widetilde
{u}\Vert_{L^{2}}^{2})\smallskip\\
&  =c_{P}^{2}\Vert\widetilde{u}\Vert_{H^{1}}^{2}\smallskip\\
\Vert P\widetilde{u}\Vert_{H^{1}} &  \leq c_{P}\Vert\widetilde{u}\Vert_{H^{1}}%
\end{align*}
and we have also $\Vert P\Vert_{H^{1}\rightarrow H^{1}}\leq c_{P}$. For
$\widetilde{u}\in\widetilde{\mathcal{V}}=P(H^{1}(B))$ we can now estimate the
minimal approximation error
\[%
\begin{array}
[c]{rcll}%
\begin{array}
[c]{c}%
\min\limits_{\widetilde{p}\in\widetilde{\mathcal{V}}_{n}}\Vert\widetilde
{u}-\widetilde{p}\Vert_{H_{1}}\\
\end{array}
& \!\!%
\begin{array}
[c]{c}%
=\\
\end{array}
& \!\!%
\begin{array}
[c]{c}%
\min\limits_{\widetilde{p}\in\widetilde{\mathcal{V}}_{n}}\Vert P\widetilde
{u}-\widetilde{p}\Vert_{H_{1}}\\
\end{array}
&
\begin{array}
[c]{l}%
P\text{ is a projection}\\
\text{and }\widetilde{u}\in\operatorname*{image}(P)
\end{array}
\smallskip\\
& \!\!= & \!\!\min\limits_{p\in\Pi_{n}}\Vert P\widetilde{u}-Pp\Vert_{H_{1}%
}\smallskip & \ \text{because }\widetilde{\mathcal{V}}_{n}=P(\Pi_{n})\\
& \!\!\leq & \!\!\min\limits_{p\in\Pi_{n}}\Vert P\Vert_{H^{1}\rightarrow
H^{1}}\Vert\widetilde{u}-p\Vert_{H_{1}}\smallskip & \\
& \!\!\leq & \!\!c_{P}\min\limits_{p\in\Pi_{n}}\Vert\widetilde{u}%
-p\Vert_{H_{1}} &
\end{array}
\]
and now we can apply the results from \cite{bbl}.

\section{Implementation \label{implementation}}

Consider the implementation of the Galerkin method of \S \ref{unique} for the
Neumann problem (\ref{en1})-(\ref{en2}) over $\Omega$ by means of the
reformulation in (\ref{en21}) over the unit ball $B_{d}$. We are to find the
function $\widetilde{u}_{n}\in\Pi_{n}$ satisfying (\ref{en21}). To do so, we
begin by selecting a basis for $\Pi_{n}$, denoting it by $\left\{  \varphi
_{1},\dots,\varphi_{N}\right\}  $, with $N\equiv N_{n}=\dim\Pi_{n}$. Generally
we use a basis that is orthonormal in the norm of $L^{2}\left(  B_{2}\right)
$. It would be better probably to use a basis that is orthonormal in the norm
of $H^{1}\left(  B_{d}\right)  $; for example, see \cite{xu2006}. \ We seek%
\begin{equation}
\widetilde{u}_{n}(x)=\sum_{k=1}^{N}\alpha_{k}\varphi_{k}(x)\label{e76}%
\end{equation}
Then (\ref{en25}) is equivalent to%
\begin{align}
&
{\displaystyle\sum\limits_{k=1}^{N}}
\alpha_{k}%
{\displaystyle\int_{B_{d}}}
\left[
{\displaystyle\sum\limits_{i,j=1}^{d}}
a_{i,j}(x)\dfrac{\partial\varphi_{k}(x)}{\partial x_{j}}\dfrac{\partial
\varphi_{\ell}(x)}{\partial x_{i}}+\gamma(x)\varphi_{k}(x)\varphi_{\ell
}(x)\right]  \left\vert \det\left[  J(x)\right]  \right\vert \,dx\medskip
\nonumber\\
&  \quad\quad\quad=%
{\displaystyle\int_{B_{d}}}
f\left(  x\right)  \varphi_{\ell}\left(  x\right)  \left\vert \det\left[
J(x)\right]  \right\vert \,dx\medskip\label{e78}\\
&  \quad\quad\quad\quad\quad\quad+%
{\displaystyle\int_{\partial B_{d}}}
g\left(  x\right)  \varphi_{\ell}\left(  x\right)  \left\vert J_{bdy}%
(x)\right\vert \,\,dx,\quad\quad\ell=1,\dots,N\nonumber
\end{align}
The function $\left\vert J_{bdy}(x)\right\vert $\ arises from the
transformation of an integral over $\partial\Omega$ to one over $\partial
B_{d}$, associated with the change from $\ell_{2}$ to $\widetilde{\ell}_{2}$
as discussed preceding (\ref{en21}). For example, in one variable the boundary
$\partial\Omega$\ \ is often represented as a mapping%
\[
\chi\left(  \theta\right)  =\left(  \chi_{1}\left(  \theta\right)  ,\chi
_{2}\left(  \theta\right)  \right)  ,\quad\quad0\leq\theta\leq2\pi.
\]
In that case, $\left\vert J_{bdy}(x)\right\vert \,$is simply $\left\vert
\chi^{\prime}\left(  \theta\right)  \right\vert $ and the associated integral
is%
\[
\int_{0}^{2\pi}g\left(  \chi\left(  \theta\right)  \right)  \varphi_{\ell
}\left(  \chi\left(  \theta\right)  \right)  \left\vert \chi^{\prime}\left(
\theta\right)  \right\vert \,d\theta
\]
In (\ref{e78}) we need to calculate the orthonormal polynomials and their
first partial derivatives; and we also need to approximate the integrals in
the linear system. For an introduction to the topic of multivariate orthogonal
polynomials, see Dunkl and Xu \cite{DX} and Xu \cite{xu2004}. For multivariate
quadrature over the unit ball in $\mathbb{R}^{d}$, see Stroud \cite{stroud}.

For the Neumann problem (\ref{en40})-(\ref{en41}) of \S \ref{nonunique}, the
implementation is basically the same. \ The basis $\left\{  \varphi_{1}%
,\dots,\varphi_{N}\right\}  $ is modified as in (\ref{en80}), with the
constant $C$ of (\ref{en81}) approximated using the quadrature in
(\ref{e106}), given below.

\subsection{The planar case}

The dimension of $\Pi_{n}$ is
\begin{equation}
N_{n}=\frac{1}{2}\left(  n+1\right)  \left(  n+2\right) \label{e79}%
\end{equation}
For notation, we replace $x$ with $\left(  x,y\right)  $. How do we choose the
orthonormal basis $\left\{  \varphi_{\ell}(x,y)\right\}  _{\ell=1}^{N}$ for
$\Pi_{n}$? Unlike the situation for the single variable case, there are many
possible orthonormal bases over $B_{d}=D$, the unit disk in $\mathbb{R}^{2}$.
We have chosen one that is particularly convenient for our computations. These
are the "ridge polynomials" introduced by Logan and\ Shepp \cite{Loga} for
solving an image reconstruction problem. We summarize here the results needed
for our work.

Let
\[
\mathcal{V}_{n}=\left\{  P\in\Pi_{n}:\left(  P,Q\right)  =0\quad\forall
Q\in\Pi_{n-1}\right\}
\]
the polynomials of degree $n$ that are orthogonal to all elements of
$\Pi_{n-1}$. Then the dimension of $\mathcal{V}_{n}$ is $n+1$; moreover,%
\begin{equation}
\Pi_{n}=\mathcal{V}_{0}\oplus\mathcal{V}_{1}\oplus\cdots\oplus\mathcal{V}%
_{n}\label{e100}%
\end{equation}
It is standard to construct orthonormal bases of each $\mathcal{V}_{n}$ and to
then combine them to form an orthonormal basis of $\Pi_{n}$ using the latter
decomposition. \ As an orthonormal basis of $\mathcal{V}_{n}$ we use%
\begin{equation}
\varphi_{n,k}(x,y)=\frac{1}{\sqrt{\pi}}U_{n}\left(  x\cos\left(  kh\right)
+y\sin\left(  kh\right)  \right)  ,\quad\left(  x,y\right)  \in D,\quad
h=\frac{\pi}{n+1}\label{e101}%
\end{equation}
for $k=0,1,\dots,n$. The function $U_{n}$ is the Chebyshev polynomial of the
second kind of degree $n$:%
\begin{equation}
U_{n}(t)=\frac{\sin\left(  n+1\right)  \theta}{\sin\theta},\quad\quad
t=\cos\theta,\quad-1\leq t\leq1,\quad n=0,1,\dots\label{e102}%
\end{equation}
The family $\left\{  \varphi_{n,k}\right\}  _{k=0}^{n}$ is an orthonormal
basis of $\mathcal{V}_{n}$. As a basis of $\Pi_{n}$, we order $\left\{
\varphi_{n,k}\right\}  $ lexicographically based on the ordering in
(\ref{e101}) and (\ref{e100}):%
\[
\left\{  \varphi_{\ell}\right\}  _{\ell=1}^{N}=\left\{  \varphi_{0,0}%
,\,\varphi_{1,0},\,\varphi_{1,1},\,\varphi_{2,0},\,\dots,\,\varphi
_{n,0},\,\dots,\varphi_{n,n}\right\}
\]
To calculate the first order partial derivatives of $\varphi_{n,k}(x,y)$, we
need $U_{n}^{^{\prime}}(t)$. The values of $U_{n}(t)$ and $U_{n}^{^{\prime}%
}(t)$ are evaluated using the standard triple recursion relations%
\begin{align*}
U_{n+1}(t)  &  =2tU_{n}(t)-U_{n-1}(t)\smallskip\\
U_{n+1}^{^{\prime}}(t)  &  =2U_{n}(t)+2tU_{n}^{^{\prime}}(t)-U_{n-1}%
^{^{\prime}}(t)
\end{align*}

For the numerical approximation of the integrals in (\ref{e78}), which are
over $B$ being the unit disk, we use the formula%
\begin{equation}
\int_{B}g(x,y)\,dx\,dy\approx\sum_{l=0}^{q}\sum_{m=0}^{2q}g\left(  r_{l}%
,\frac{2\pi\,m}{2q+1}\right)  \omega_{l}\frac{2\pi}{2q+1}r_{l}\label{e106}%
\end{equation}
Here the numbers $\omega_{l}$ are the weights of the $\left(  q+1\right)
$-point Gauss-Legendre quadrature formula on $[0,1]$. Note that
\[
\int_{0}^{1}p(x)dx=\sum_{l=0}^{q}p(r_{l})\omega_{l},
\]
for all single-variable polynomials $p(x)$ with $\deg\left(  p\right)
\leq2q+1 $. The formula (\ref{e106}) uses the trapezoidal rule with $2q+1$
subdivisions for the integration over $\overline{B}_{d}$ in the azimuthal
variable. This quadrature is exact for all polynomials $g\in\Pi_{2q}$. This
formula is also the basis of the hyperinterpolation formula discussed in
\cite{hac}.

\subsection{The three dimensional case}

In the three dimensional case the dimension of $\Pi_{n}$ is given by
\[
N_{n}={\binom{n+3}{3}}%
\]
and we choose the following orthogonal polynomials on the unit ball
\begin{align}
\varphi_{m,j,\beta}(x) &  =c_{m,j}p_{j}^{(0,m-2j+\frac{1}{2})}(2\Vert
x\Vert^{2}-1)S_{\beta,m-2j}\left(  x\right)  \medskip\nonumber\\
&  =c_{m,j}\Vert x\Vert^{m-2j}p_{j}^{(0,m-2j+\frac{1}{2})}(2\Vert x\Vert
^{2}-1)S_{\beta,m-2j}\left(  \frac{x}{\Vert x\Vert}\right)  ,\medskip
\label{eq1000}\\
j &  =0,\ldots,\lfloor m/2\rfloor,\quad\beta=0,1,\ldots,2(m-2j),\quad
m=0,1,\ldots,n\nonumber
\end{align}
The constants $c_{m,j}$ are given by $c_{m,j}=2^{\frac{5}{4}+\frac{m}{2}-j}$;
and the functions $p_{j}^{(0,m-2j+\frac{1}{2})}$ are the normalized Jacobi
polynomials. The functions $S_{\beta,m-2j}$ are spherical harmonic functions
and they are orthonormal on the sphere $\mathbb{S}^{2}\subset\mathbb{R}^{3}$.
See \cite{DX,ACH-Dir} for the definition of these functions. In \cite{ACH-Dir}
one also finds the quadrature methods which we use to approximate the
integrals over $B_{1}(0)$ in (\ref{en20}) and (\ref{en21}). The functional
$\widetilde{\ell}_{2}$ in (\ref{en21}) is given by%
\begin{align}
\widetilde{\ell}_{2}(v)  & =%
{\displaystyle\int_{0}^{\pi}}
{\displaystyle\int_{0}^{2\pi}}
g(\Phi(\Upsilon(1,\theta,\phi)))\label{eq1001}\\
& \cdot\Vert(\Phi\circ\Upsilon)_{\theta}(1,\theta,\phi)\times(\Phi
\circ\Upsilon)_{\phi}(1,\theta,\phi)\Vert\,v(\Phi(\Upsilon(1,\theta
,\phi)))d\phi\,d\theta\nonumber
\end{align}
where
\begin{equation}
\Upsilon(\rho,\theta,\phi):=\rho\left(  \sin(\theta)\cos(\phi),\sin
(\theta)\sin(\phi),\cos(\theta)\right) \label{eq1002}%
\end{equation}
is the usual transformation between spherical and Cartesian coordinates and
the indices denote the partial derivatives. For the numerical approximation of
the integral in (\ref{eq1001}) we use traezoidal rules in the $\phi$ direction
and Gau{\ss }-Legendre formulas for the $\theta$ direction.

\section{Numerical examples \label{examples}}

The construction of our examples is very similar to that given in
\cite{ACH-Dir} for the Dirichlet problem. Our first two transformations $\Phi$
have been so chosen that we can invert explicitly the mapping $\Phi$, to be
able to better construct our test examples. \textbf{This is not needed when
applying the method}; but it simplifies the construction of our test cases.
Given $\Phi$, we need to calculate analytically the matrix%
\begin{equation}
A(x)=J\left(  x\right)  ^{-1}J\left(  x\right)  ^{-\text{T}}.\label{e131}%
\end{equation}

\subsection{The planar case}%

\begin{figure}
[tb]
\begin{center}
\includegraphics[
height=3in,
width=3.9998in
]%
{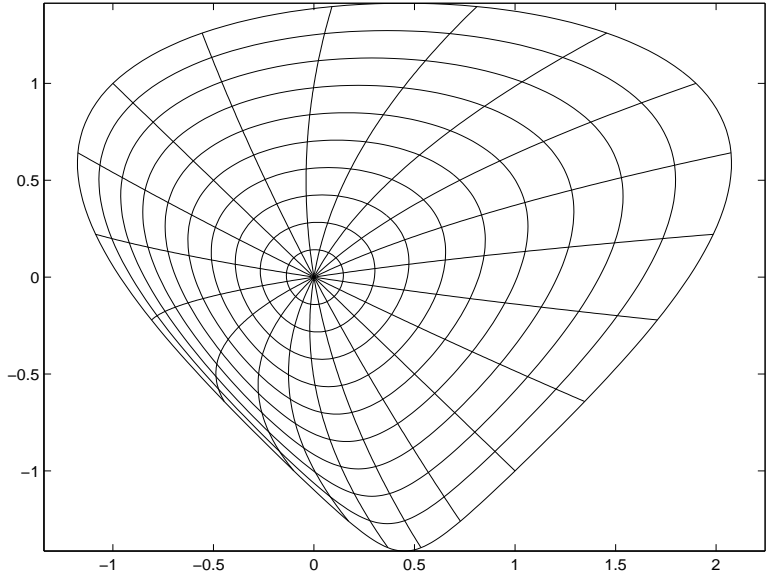}%
\caption{Images of (\ref{e132}), with $a=0.5$, for lines of constant radius
and constant azimuth on the unit disk.}%
\label{fig1}%
\end{center}
\end{figure}
For our variables, we replace a point $x\in B_{d}$ with $\left(  x,y\right)
$, and we replace a point $s\in\Omega$ with $\left(  s,t\right)  $. \ Define
the mapping $\Phi:\overline{B}\rightarrow\overline{\Omega}$ by $\left(
s,t\right)  =\Phi\left(  x,y\right)  $,%
\begin{equation}%
\begin{array}
[c]{l}%
s=x-y+ax^{2}\\
t=x+y
\end{array}
\label{e132}%
\end{equation}
with $0<a<1$. It can be shown that $\Phi$ is a 1-1 mapping from the unit disk
$\overline{B}$. In particular, the inverse mapping $\Psi:\overline{\Omega
}\rightarrow\overline{B}$ is given by%
\begin{equation}%
\begin{array}
[c]{l}%
x=\dfrac{1}{a}\left[  -1+\sqrt{1+a\left(  s+t\right)  }\right]  \medskip\\
y=\dfrac{1}{a}\left[  at-\left(  -1+\sqrt{1+a\left(  s+t\right)  }\right)
\right]
\end{array}
\label{e133}%
\end{equation}
In Figure \ref{fig1}, we give the images in $\overline{\Omega}$ of the circles
$r=j/10$, $j=1,\dots,10$ and the azimuthal lines $\theta=j\pi/10$,
$j=1,\dots,20$.

The following information is needed when implementing the transformation from
$-\Delta u+\gamma u=f$ on $\Omega$ to a new equation on $B$:%
\[
D\Phi=J\left(  x,y\right)  =\left(
\begin{array}
[c]{cc}%
1+2ax & -1\\
1 & 1
\end{array}
\right)
\]%
\[
\det\left(  J\right)  =2\left(  1+ax\right)
\]%
\[
J\left(  x\right)  ^{-1}=\frac{1}{2\left(  1+ax\right)  }\left(
\begin{array}
[c]{cc}%
1 & 1\\
-1 & 1+2ax
\end{array}
\right)
\]%
\[
A=J\left(  x\right)  ^{-1}J\left(  x\right)  ^{-\text{T}}=\frac{1}{2\left(
1+ax\right)  ^{2}}\left(
\begin{array}
[c]{cc}%
1 & ax\\
ax & 2a^{2}x^{2}+2ax+1
\end{array}
\right)
\]
The latter are the coefficients needed to define $\widetilde{\mathcal{A}}$ in
(\ref{en20}).%

\begin{figure}
[tb]
\begin{center}
\includegraphics[
height=3in,
width=3.9998in
]%
{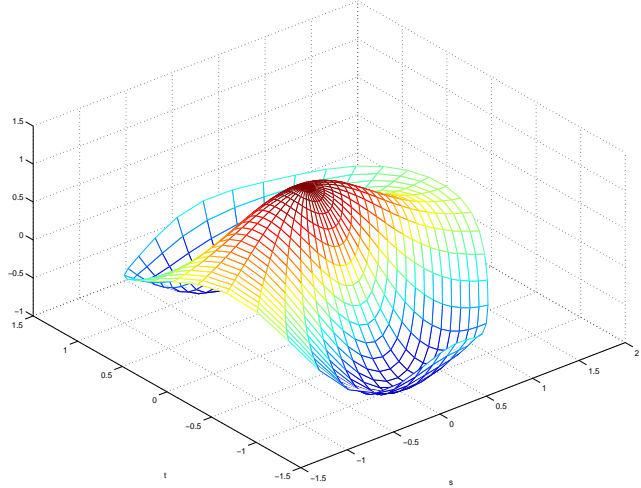}%
\caption{The function $u(s,t)$ of (\ref{e138})}%
\label{true_soln}%
\end{center}
\end{figure}

We give numerical results for solving the equation%
\begin{equation}
-\Delta u\left(  s,t\right)  +e^{s-t}u\left(  s,t\right)  =f\left(
s,t\right)  ,\quad\quad\left(  s,t\right)  \in\Omega\label{e136}%
\end{equation}
As a test case, we choose%
\begin{equation}
u\left(  s,t\right)  =e^{-s^{2}}\cos\left(  \pi t\right)  ,\quad\quad\left(
s,t\right)  \in\Omega\label{e138}%
\end{equation}
The solution is pictured in Figure \ref{true_soln}. To find $f(s,t)$, we use
(\ref{e136}) and (\ref{e138}). We use the domain parameter $a=0.5$, with
$\Omega$ pictured in Figure \ref{fig1}.

\begin{figure}
[tb]
\begin{center}
\includegraphics[
height=3in,
width=3.9998in
]%
{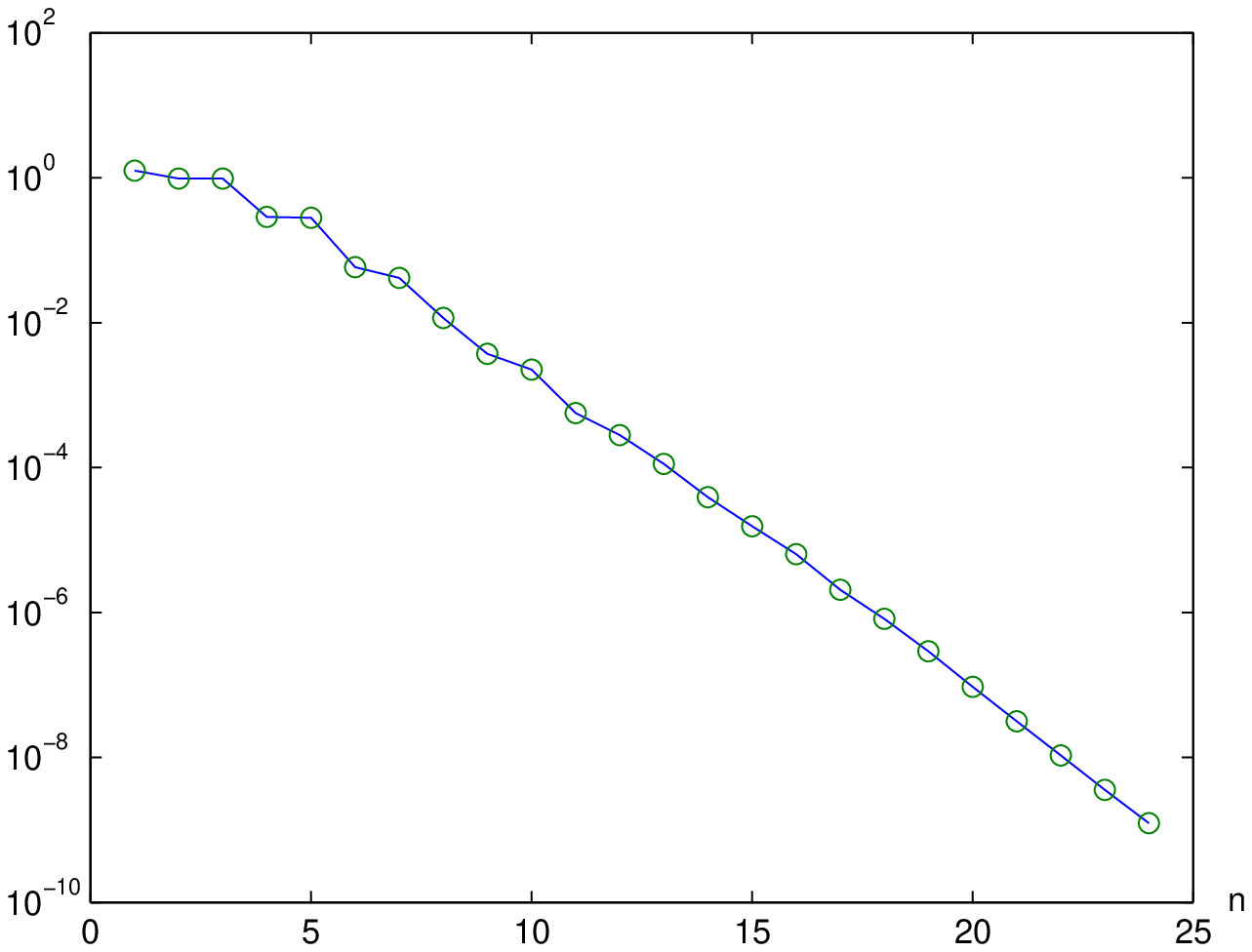}%
\caption{Errors from Table \ref{table1}}%
\label{error_planar}%
\end{center}
\end{figure}

Numerical results are given in Table \ref{table1} for even values of $n$.
\ The integrations in (\ref{e78}) were performed with (\ref{e106}); and the
integration parameter $q$ ranged from $10$ to $30$. We give the condition
numbers of the linear system (\ref{e78}) as produced in \textsc{Matlab}. To
calculate the error, we evaluate the numerical solution and the error on the
grid
\begin{align*}
\Phi\left(  x_{i,j},y_{i,j}\right)   &  =\Phi\left(  r_{i}\cos\theta_{j}%
,r_{i}\sin\theta_{j}\right)  \smallskip\\
\left(  r_{i},\theta_{j}\right)   &  =\left(  \frac{i}{10},\frac{j\pi}%
{10}\right)  ,\quad\quad i=0,1,\dots10;\quad j=1,\dots20
\end{align*}
The results are shown graphically in Figure \ref{error_planar}. The use of a
semi-log scale demonstrates the exponential convergence of the method as the
degree increases.%

\begin{table}[tb] \centering
\caption{Maximum errors in Galerkin solution $u_n$}\label{table1}%
\begin{tabular}
[c]{|c|c|c|c||c|c|c|c|}\hline
$n$ & $N_{n}$ & $\left\Vert u-u_{n}\right\Vert _{\infty}$ & \textit{cond} &
$n$ & $N_{n}$ & $\left\Vert u-u_{n}\right\Vert _{\infty}$ & \textit{cond}%
\\\hline
$2$ & $6$ & $9.71E-1$ & $14.5$ & $14$ & $120$ & $3.90E-5$ & $6227$\\\hline
$4$ & $15$ & $2.87E-1$ & $86.1$ & $16$ & $153$ & $6.37E-6$ & $10250$\\\hline
$6$ & $28$ & $5.85E-2$ & $309$ & $18$ & $190$ & $8.20E-7$ & $15960$\\\hline
$8$ & $45$ & $1.16E-2$ & $824$ & $20$ & $231$ & $9.44E-8$ & $23770$\\\hline
$10$ & $66$ & $2.26E-3$ & $1819$ & $22$ & $276$ & $1.06E-8$ & $34170$\\\hline
$12$ & $91$ & $2.81E-4$ & $3527$ & $24$ & $325$ & $1.24E-9$ & $47650$\\\hline
\end{tabular}%
\end{table}%
%

%

\begin{figure}
[tb]
\begin{center}
\includegraphics[
height=3in,
width=3.9998in
]%
{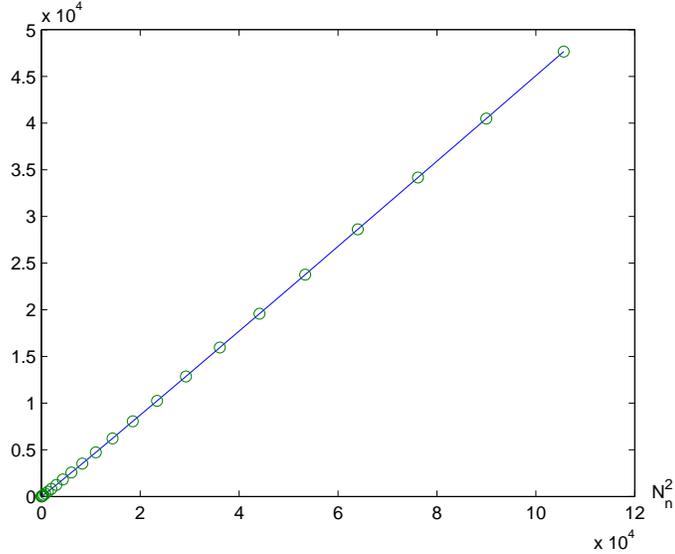}%
\caption{Condition numbers from Table \ref{table1}}%
\label{cond_planar}%
\end{center}
\end{figure}

To examine experimentally the behaviour of the condition numbers for the
linear system (\ref{e78}), we have graphed the condition numbers from Table
\ref{table1} in Figure \ref{cond_planar}. Note that we are graphing $N_{n}%
^{2}$ vs. the condition number of the associated linear system. \ The graph
seems to indicate that the condition number of the system (\ref{e78}) is
directly proportional to the square of the order of the system, with the order
given in (\ref{e79}).

For the Poisson equation%
\[
-\Delta u\left(  s,t\right)  =f\left(  s,t\right)  ,\quad\quad\left(
s,t\right)  \in\Omega
\]
with the same true solution as in (\ref{e138}), we use the numerical method
given in \S \ref{nonunique}. The numerical results are comparable. \ For
example, with $n=20$, we obtain $\left\Vert u-u_{n}\right\Vert _{\infty
}=9.90\times10^{-8}$ and the condition number is approximately $14980$.%
\begin{figure}
[tb]
\begin{center}
\includegraphics[
height=3in,
width=3.9998in
]%
{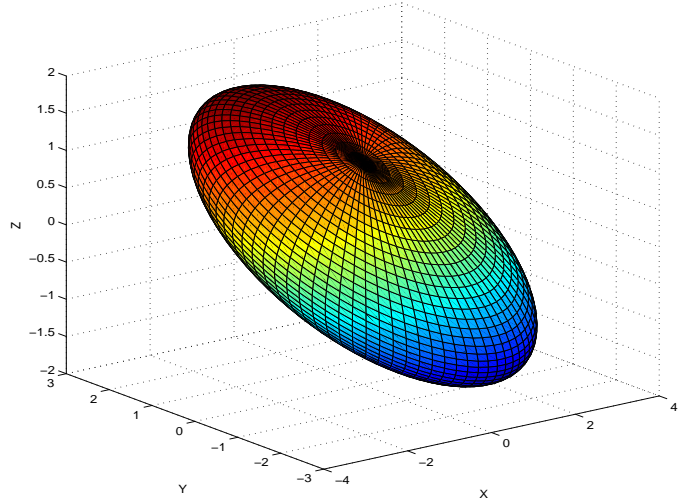}%
\caption{The boundary of $\Omega_{1}$}%
\label{figure_ellipsoid}%
\end{center}
\end{figure}

\subsection{The three dimensional case}

To illustrate that the proposed spectral method converges rapidly, we first
use a simple test example. We choose the linear transformation
\[
s:=\Phi_{1}(x)=\left(
\begin{array}
[c]{c}%
x_{1}-3x_{2}\\
2x_{1}+x_{2}\\
x_{1}+x_{2}+x_{3}%
\end{array}
\right)  ,
\]
so that $B_{1}(0)$ is transformed to an ellipsoid $\Omega_{1}$; see figure
\ref{figure_ellipsoid}. For this transformation $D\Phi_{1}$ and $J_{1}%
=\det(D\Phi_{1})$ are constant functions. For a test solution, we use the
function
\begin{equation}
u(s)=s_{1}e^{s_{2}}\sin(s_{3})\label{eq2000}%
\end{equation}
which is analytic in each variable.%

\begin{table}[tb] \centering
\caption{Maximum errors in Galerkin solution $u_n$}\label{table2}%
\begin{tabular}
[c]{|c|c|c|c||c|c|c|c|}\hline
$n$ & $N_{n}$ & $\left\Vert u-u_{n}\right\Vert _{\infty}$ & \textit{cond} &
$n$ & $N_{n}$ & $\left\Vert u-u_{n}\right\Vert _{\infty}$ & \textit{cond}%
\\\hline
$1$ & $4$ & $9.22E+00$ & $8$ & $9$ & $220$ & $4.15E-04$ & $1964$\\\hline
$2$ & $10$ & $5.25E+00$ & $31$ & $10$ & $286$ & $6.84E-05$ & $2794$\\\hline
$3$ & $20$ & $1.92E+00$ & $79$ & $11$ & $364$ & $1.11E-05$ & $3862$\\\hline
$4$ & $35$ & $5.80E-01$ & $167$ & $12$ & $455$ & $1.60E-06$ & $5211$\\\hline
$5$ & $56$ & $1.62E-01$ & $314$ & $13$ & $560$ & $2.06E-07$ & $6888$\\\hline
$6$ & $84$ & $4.53E-02$ & $540$ & $14$ & $680$ & $2.60E-08$ & $8937$\\\hline
$7$ & $120$ & $1.03E-02$ & $871$ & $15$ & $816$ & $3.01E-09$ & $11415$\\\hline
$8$ & $165$ & $2.31E-03$ & $1335$ & $16$ & $969$ & $3.13E-10$ & $14376$%
\\\hline
\end{tabular}%
\end{table}%
%

\begin{figure}
[tb]
\begin{center}
\includegraphics[
height=3in,
width=3.9998in
]%
{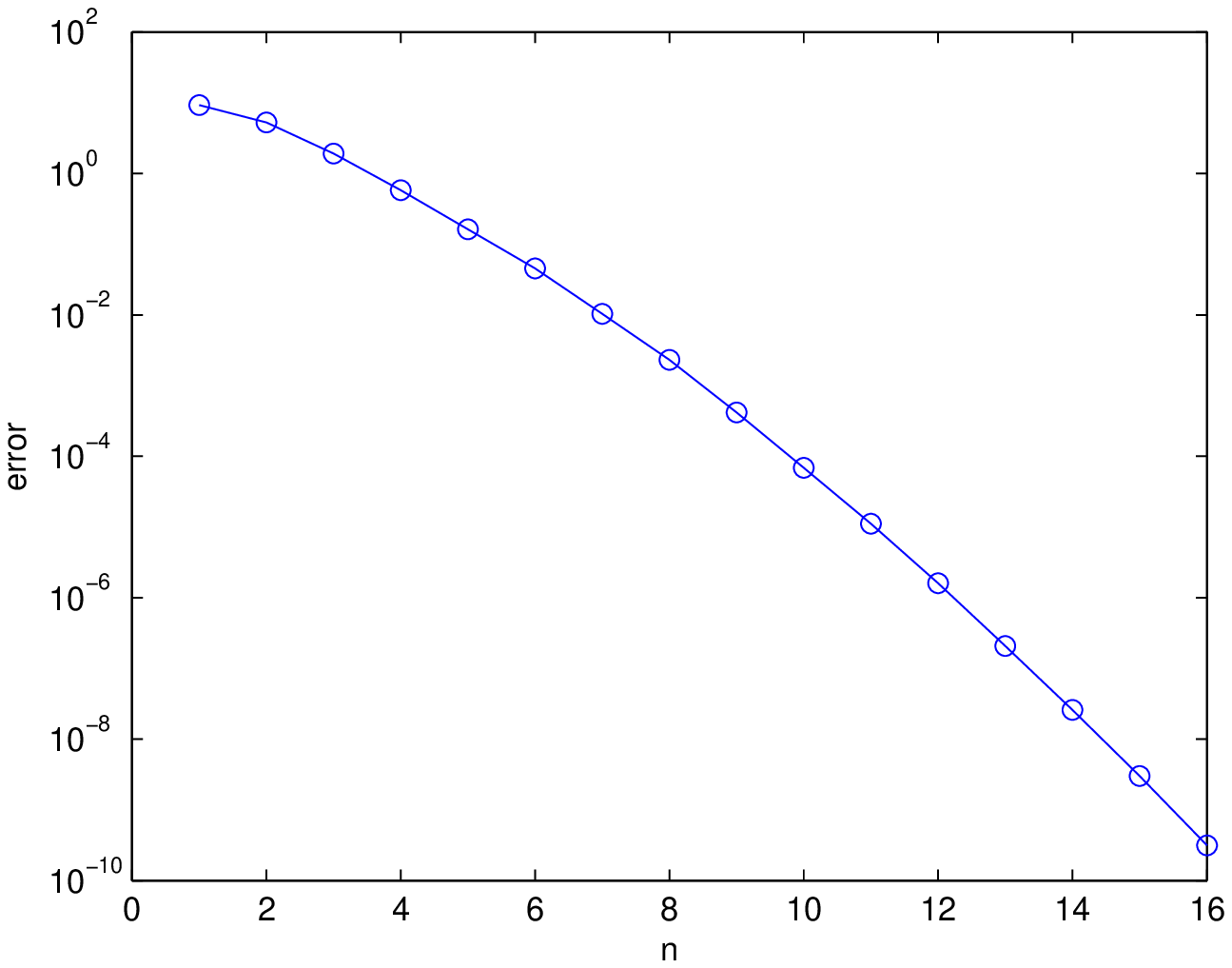}%
\caption{Errors from Table \ref{table2}}%
\label{figure_omega1e}%
\end{center}
\end{figure}
\begin{figure}
[tb]
\begin{center}
\includegraphics[
height=3in,
width=3.9998in
]%
{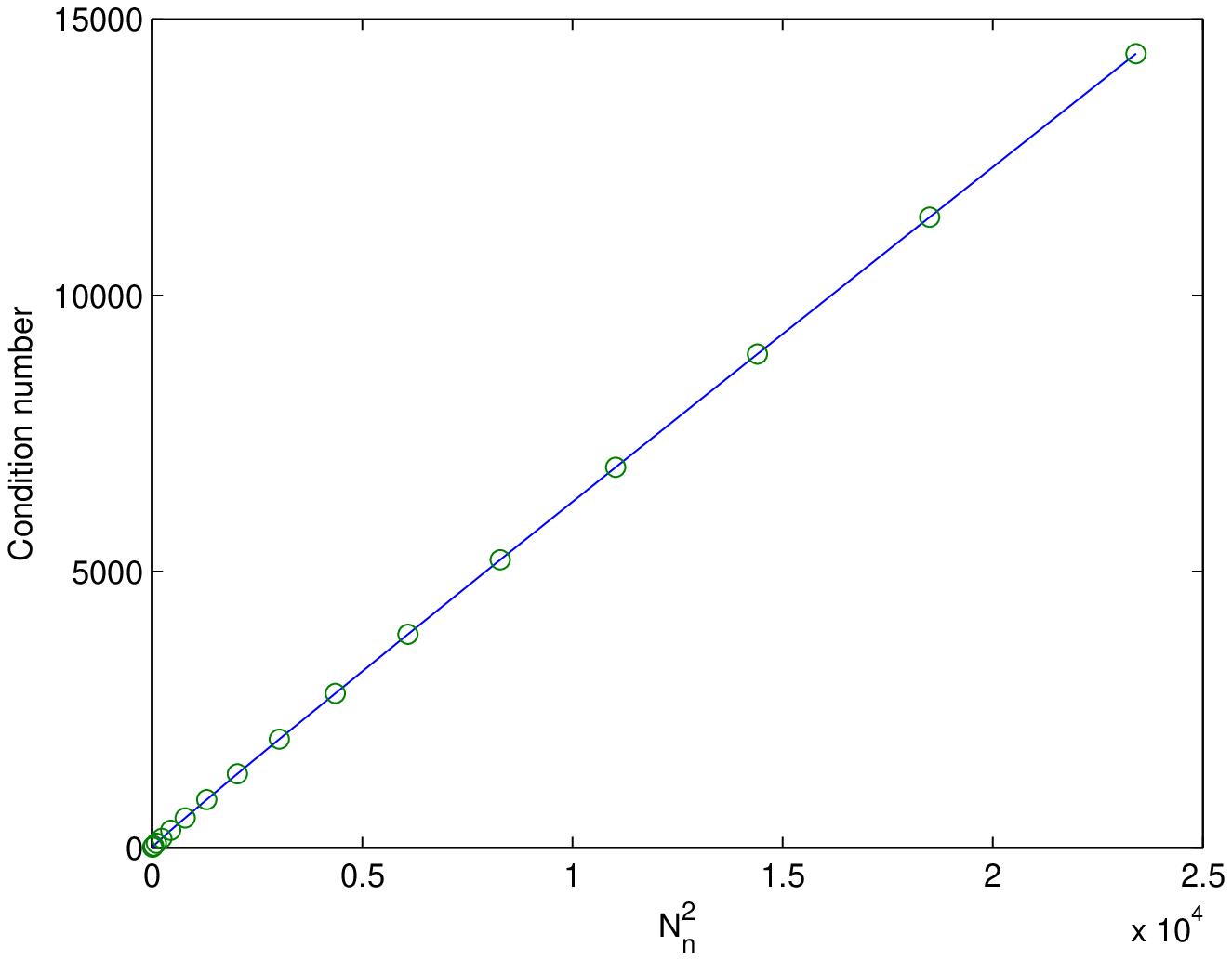}%
\caption{Conditions numbers from Table \ref{table2}}%
\label{figure_omega1c}%
\end{center}
\end{figure}

Table \ref{table2} shows the errors and the development of the condition
numbers for the solution of (\ref{en1}) on $\Omega_{1}$. The associated graphs
for the errors and condition numbers are shown in figures \ref{figure_omega1e}
and \ref{figure_omega1c}, respectively. The graph of the error is consistent
with exponential convergence; and the condition number seems to have a growth
proportional to the square of the number of degrees of freedom $N_{n}$.

Next we study domains $\Omega$ which are star shaped with respect to the
origin,
\begin{equation}
\Omega_{2}=\{x\in\mathbb{R}^{3}\mid x=\Upsilon(\rho,\theta,\phi),\quad
0\leq\rho\leq R(\theta,\phi)\}.\label{eq2001a}%
\end{equation}
See (\ref{eq1002}) for the definition of $\Upsilon$, and $R:\mathbb{S}%
^{2}\rightarrow(0,\infty)$ is assumed to be a $C^{\infty}$ function. In this
case we can construct arbitrarily smooth and invertible mappings $\Phi
:B_{1}(0)\rightarrow\Omega_{2}$ as we will show now. First we define a
function $t:[0,1]\rightarrow\lbrack0,1]$
\begin{equation}
t(\rho):=\left\{
\begin{array}
[c]{cc}%
0, & 0\leq\rho\leq\frac{1}{2},\\
2^{e_{s}}(\rho-\frac{1}{2})^{e_{s}}, & \frac{1}{2}<\rho\leq1.
\end{array}
\right. \label{eq2002}%
\end{equation}
the parameter $e_{s}\in\mathbb{N}$ determines the smoothness of $t\in
C^{e_{s}-1}[0,1]$. For the following we will assume that $R(\theta,\phi)>1$,
for all $\theta$ and $\phi$; this follows after an appropriate scaling of the
problem. With the help of $t$ we define the function $\widetilde{R}$ which is
monotone increasing from $0$ to $R(\theta,\phi)$ on $[0,1]$ and equal to the
identity on $[0,1/2]$,
\[
\widetilde{R}(\rho,\theta,\phi):=t(\rho)R(\theta,\phi)+(1-t(\rho))\rho
\]
Because
\[
\frac{\partial}{\partial\rho}\widetilde{R}(\rho,\theta,\phi)=t^{\prime}%
(\rho)(R(\theta,\phi)-\rho)+(1-t(\rho))>0,\quad\rho\in\lbrack0,1]
\]
the function $\widetilde{R}$ is an invertible function of $\rho$ of class
$C^{e_{s}-1}$. The transformation $\Phi_{2}:B_{1}(0)\rightarrow\Omega_{2}$ is
defined by
\[
\Phi_{2}(x):=\Upsilon(\widetilde{R}(\rho,\theta,\phi),\theta,\phi),\quad
x=\Upsilon(\rho,\theta,\phi)\in B_{1}(0)
\]
The properties of $\widetilde{R}$ imply that $\Phi_{2}$ is equal to the
identity on $B_{\frac{1}{2}}(0)$ and the outside shell $B_{1}(0)\setminus
B_{\frac{1}{2}}(0)$ is deformed by $\Phi_{2}$ to cover $\Omega_{2}\setminus
B_{\frac{1}{2}}(0)$.%

\begin{figure}
[tb]
\begin{center}
\includegraphics[
height=3in,
width=3.9998in
]%
{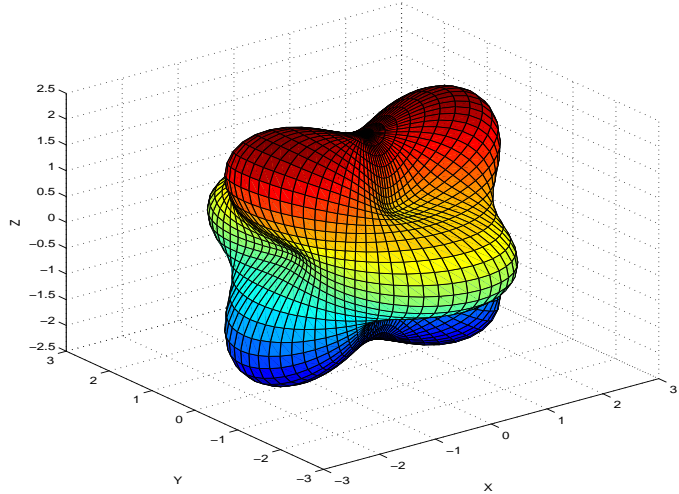}%
\caption{A view of $\partial\Omega_{2}$}%
\label{shape_omega2y}%
\end{center}
\end{figure}
%

\begin{figure}
[tb]
\begin{center}
\includegraphics[
height=3in,
width=3.9998in
]%
{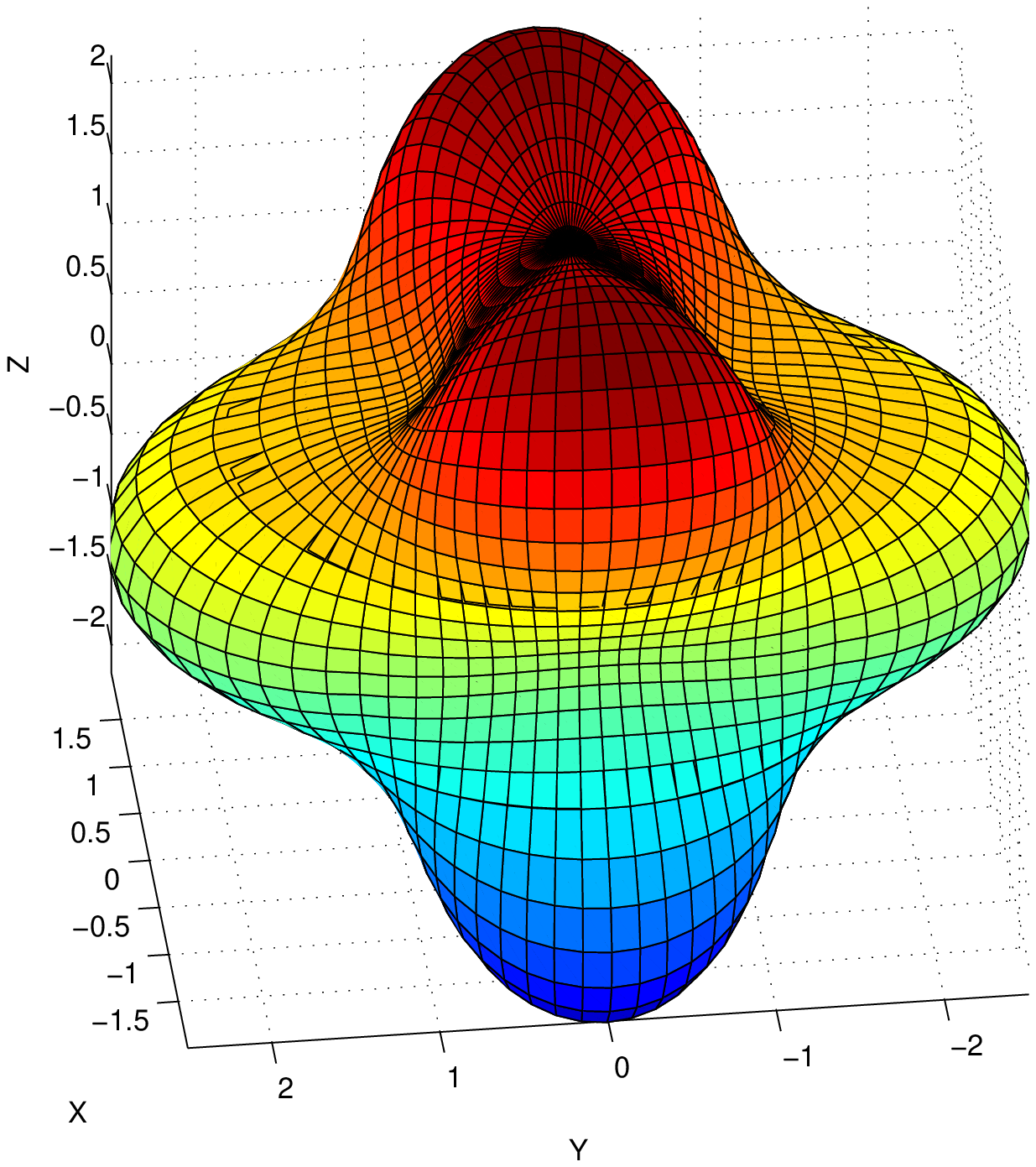}%
\caption{Another view of $\partial\Omega_{2}$}%
\label{shape_omega2z}%
\end{center}
\end{figure}

For a test surface, we use
\begin{align}
R(\theta,\phi)  & =2+\frac{3}{4}\cos(2\phi)\sin(\theta)^{2}(7\cos(\theta
)^{2}-1)\label{eq2002a}\\
e_{s}  & =5;\nonumber
\end{align}
see figures \ref{shape_omega2y}-\ref{shape_omega2z} for pictures of
$\partial\Omega_{2}$. For our test example, we use $u$ from (\ref{eq2000}).

The term $\cos(2\phi)\sin(\theta)^{2}(7\cos(\theta)^{2}-1)$ is a spherical
harmonic function which shows $R\in C^{\infty}(\mathbb{S}^{2})$, and the
factor $3/4$ is used to guarantee $R>1$. For the transformation $\Phi_{2}$ we
get $\Phi_{2}\in C^{4}(B_{1}(0))$, so we expect a convergence of order
$O(n^{-4})$. Our spectral method will now approximate $u\circ\Phi_{2}$ on the
unit ball, which varies much more than the function in our first example.

We also note that one might ask why we do not further increase $e_{s}$ (see
(\ref{eq2002})) to get a better order of convergence. It is possible to do
this, but the price one pays is in larger derivatives of $u\circ\Phi_{2}$, and
this may result in larger errors for the range of $n$ values where we actually
calculate the approximation. The search for an optimal $e_{s}$ is a problem on
its own, but it also depends on the solution $u$. So we have chosen $e_{s}=5$
in order to demonstrate our method, showing that the qualitative behaviour of
the error is the same as in our earlier examples.

\begin{figure}
[tb]
\begin{center}
\includegraphics[
height=3in,
width=3.9998in
]%
{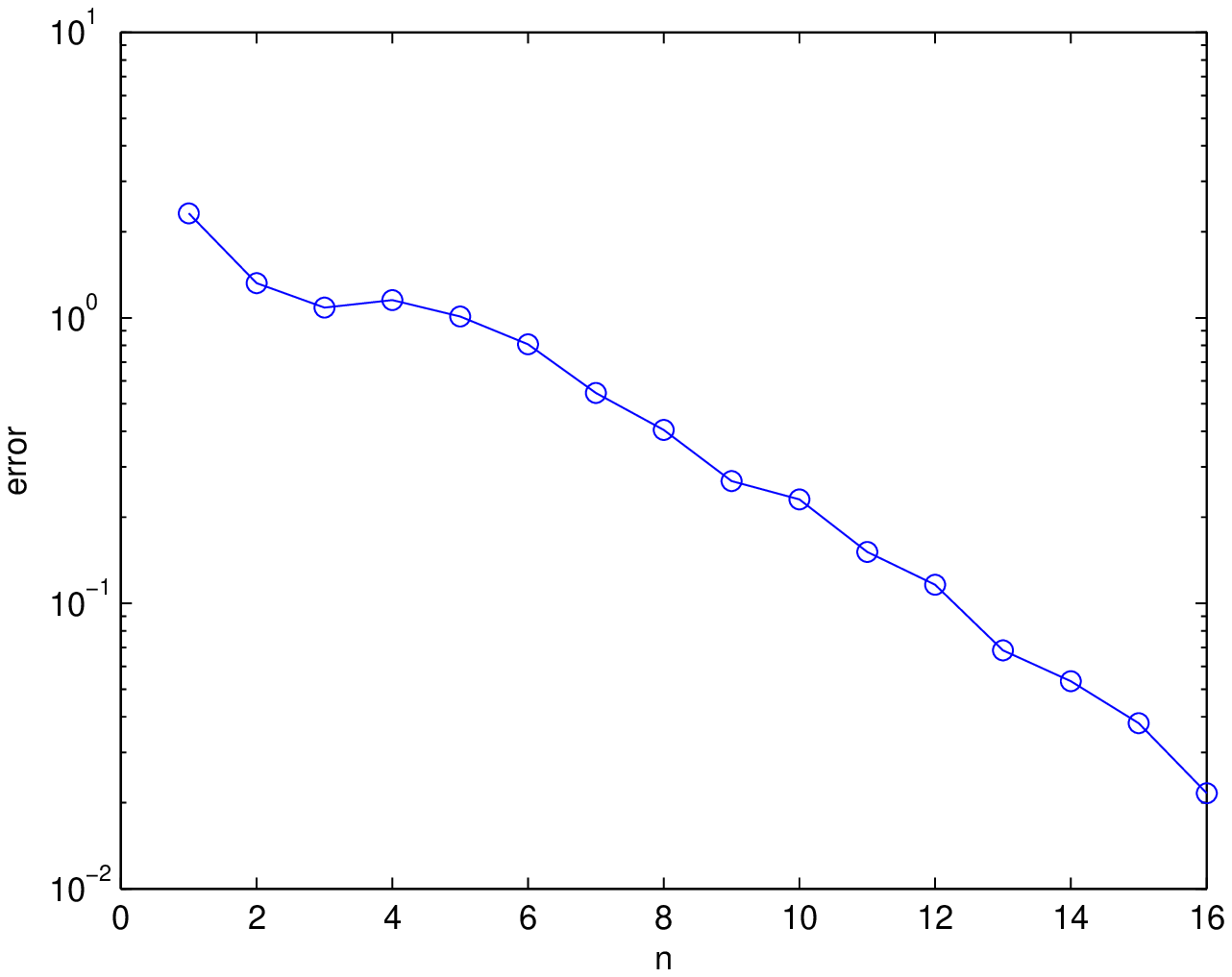}%
\caption{Errors from table \ref{table3}}%
\label{figure_omega2e}%
\end{center}
\end{figure}

\begin{figure}
[tb]
\begin{center}
\includegraphics[
height=3in,
width=3.9998in
]%
{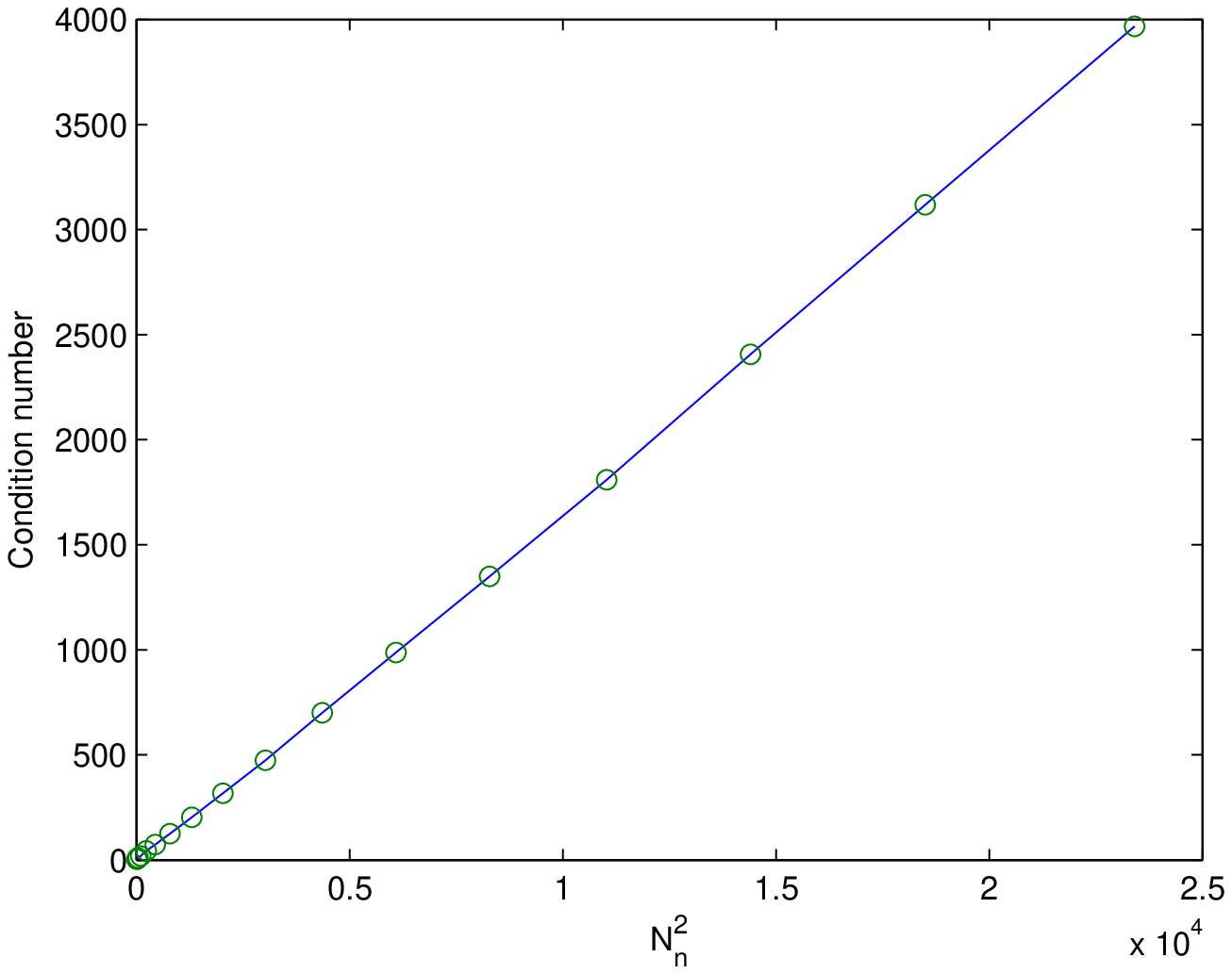}%
\caption{Condition numbers from table \ref{table3}}%
\label{figure_omega2c}%
\end{center}
\end{figure}

The results of our calculations are given in table \ref{table3}, and the
associated graphs of the errors and condition numbers are shown in figures
\ref{figure_omega2e} and \ref{figure_omega2c}, respectively. The graph in
Figure \ref{figure_omega2c} shows that the condition numbers of the systems
grow more slowly than in our first example, but again the condition numbers
appear to be proportional to $N_{n}^{2}$. The graph of the error in Figure
\ref{figure_omega2e} again resembles a line and this implies exponential
convergence; but the line has a much smaller slope than in the first example
so that the error is only reduced to about $0.02$ when we use degree $16$.
What we expect is a convergence of order $O(n^{-4})$, but the graph does not
reveal this behavior in the range of $n$ values we have used. Rather, the
convergence appears to be exponential. In the future we plan on repeating this
numerical example with an improved extension $\Phi$ of the boundary given in
(\ref{eq2002a}).

When given a mapping $\varphi:\partial B\rightarrow\partial\Omega$, it is
often nontrivial to find an extension $\Phi:\overline{B}\underset
{onto}{\overset{1-1}{\longrightarrow}}\overline{\Omega}$ with $\left.
\Phi\right\vert _{\partial B}=\varphi$ and with other needed properties. For
example, consider a star-like region $\Omega$ whose boundary surface
$\partial\Omega$ is given by%
\[
\rho=R(\theta,\phi)
\]
with $R:\mathbb{S}^{2}\rightarrow\partial\Omega$. It might seem natural to use%
\[
\Phi\left(  \rho,\theta,\phi\right)  =\rho R(\theta,\phi),\quad0\leq\rho
\leq1,\quad0\leq\theta\leq\pi,\quad0\leq\phi\leq2\pi
\]
However, such a function $\Phi$ is not continuously differentiable at $\rho
=0$. We are exploring this general problem, looking at ways of producing
$\Phi$ with the properties that are needed for implementing our spectral
method.\bigskip

\noindent\textsc{ACKNOWLEDGEMENTS}. The authors would like to thank Professor
Weimin Han for his careful proofreading of the manuscript.%

\begin{table}[tb] \centering
\caption{Maximum errors in Galerkin solution $u_n$}\label{table3}%
\begin{tabular}
[c]{|c|c|c|c||c|c|c|c|}\hline
$n$ & $N_{n}$ & $\left\Vert u-u_{n}\right\Vert _{\infty}$ & \textit{cond} &
$n$ & $N_{n}$ & $\left\Vert u-u_{n}\right\Vert _{\infty}$ & \textit{cond}%
\\\hline
$1$ & $4$ & $2.322$ & $3$ & $9$ & $220$ & $0.268$ & $475$\\\hline
$2$ & $10$ & $1.321$ & $10$ & $10$ & $286$ & $0.231$ & $701$\\\hline
$3$ & $20$ & $1.085$ & $19$ & $11$ & $364$ & $0.151$ & $987$\\\hline
$4$ & $35$ & $1.152$ & $44$ & $12$ & $455$ & $0.116$ & $1350$\\\hline
$5$ & $56$ & $1.010$ & $73$ & $13$ & $560$ & $0.068$ & $1809$\\\hline
$6$ & $84$ & $0.807$ & $125$ & $14$ & $680$ & $0.053$ & $2406$\\\hline
$7$ & $120$ & $0.545$ & $203$ & $15$ & $816$ & $0.038$ & $3118$\\\hline
$8$ & $165$ & $0.404$ & $318$ & $16$ & $969$ & $0.022$ & $3967$\\\hline
\end{tabular}%
\end{table}%
%

%

\end{document}